\documentclass{article}
\usepackage{tikz}
\usepackage{booktabs}
\usepackage{caption}
\usepackage{subcaption}
\usepackage{pgffor}
\usepackage{multicol}
\usepackage{graphicx}
\usepackage{arxiv}
\usepackage{algorithm}
\usepackage{array} 
\usepackage{amsmath}
\usepackage{algpseudocode}
\usepackage[utf8]{inputenc} 
\usepackage[T1]{fontenc}    
\usepackage{hyperref}       
\usepackage{url}            
\usepackage{booktabs}       
\usepackage{amsfonts}       
\usepackage{nicefrac}       
\usepackage{microtype}      
\usepackage{lipsum}
\usepackage{graphicx}
\graphicspath{ {./images/} }

\title{Multi-scale representation of integer sets: application to prime numbers}

\author{
 Mahmoud Melkemi \\
  Université de Haute Alsace \\
  France\\
  \texttt{mahmoud.melkemi@uha.fr} \\
}

\begin{document}

\maketitle
\begin{abstract}

Consider the partition of the set of natural numbers $\mathbb{N}$ into consecutive blocks, each containing $8^k$ integers. 
Denote these blocks by $P^{(k)}(j)$, where $j \in \mathbb{N}$.

Partition each block $P^{(k)}(j)$ into eight disjoint subsets, denoted $I_1, I_2, \ldots, I_8$, each consisting of $8^{k-1}$ consecutive integers. Define
\begin{equation*}
\phi(I_i)=
\begin{cases}
1 & \text{if } I_i \text{ contains at least one prime number},\\
0 & \text{otherwise}.
\end{cases}
\label{eq:phiMoins1}
\end{equation*}

The integer $n=(\phi(I_1)\phi(I_2)\cdots\phi(I_8))_2$
is therefore a number between 0 and 255. It encodes the configuration indicating whether each subset $I_i$ contains at least one prime 
number.
Among the set of integers less than or equal to $m$, we study the quantity $\pi_m^{(k)}(n)$, defined as the number of blocks $P^{(k)}(j)$ 
having configuration $n$. More precisely,
\begin{equation*}
\pi^{(k)}_{m}(n)=
\operatorname{card}
\left\{
P^{(k)}(j) \mid
(\phi(I_1)\phi(I_2)\cdots\phi(I_8))_2=n,\;
j=0,\ldots,\frac{m}{8^k}
\right\}.
\end{equation*}

This article provides estimates of these quantities for every integer $k$. 
To achieve this goal, we introduce a multi-scale representation of the integers for studying arithmetic properties of sets of integers. Information is organized as a hierarchy of nested sequences, where each level highlights a particular characteristic of the property under investigation. Here, we apply this framework to the distribution of prime numbers to investigate the problem described above.

Although this work does not claim any new breakthrough on this classical problem, the proposed multi-scale representation reveals structural features of the distribution of prime numbers. Moreover, the natural numbers are mapped into hierarchical sequences taking values in $\{0,\ldots,255\}$, providing a compact representation that facilitates the analysis of the encoded property across different scales. 
Instead of manipulating large blocks of integers directly, we encode them using small integers between 0 and 255, making the underlying patterns easier to visualize and interpret.

\end{abstract}


\section{Introduction}
This article presents a multi-scale representation of the natural numbers that enables the analysis of their arithmetic properties at different resolutions. This framework provides a unified approach for simultaneously studying local structures (fine-scale properties of integers) and global structures (their overall behavior).

Multi-scale methods have proved effective in many scientific fields, particularly in signal and image processing. In number theory, multi-scale ideas have also been applied to prime numbers~\cite{Iovane2009}, with connections to fractal geometry and numerical accelerators via the Golden Mean. Our approach differs from these existing methods and is designed for the study of general arithmetic properties of integers.

To illustrate the proposed framework, we focus on primality, a classical topic that has seen remarkable theoretical advances~\cite{dusart2018,goldston2009,maynard2015}. 
Our objective is not to address these deep theoretical results but rather to demonstrate how a multi-scale representation can reveal arithmetic properties through pattern analysis across different scales.

The proposed representation generates, at each scale, a sequence of integers encoding the property under study.
In this article, every encoded value belongs to the interval $[0,255]$.

Each integer compactly summarizes the distribution of the property under consideration over a large block of consecutive integers.

To illustrate the construction, we consider the property of primality on $\mathbb{N}^*$.

At scale $1$, $\mathbb{N}^*$ is viewed as the union of the sets
\[
P^{(1)}_n = \{8n+1,\ldots,8n+8\}, \qquad n\in\mathbb{N}.
\]

More generally, for every scale $k\ge2$, the sets are defined recursively by

\begin{equation}
P^{(k)}_n=
\bigcup_{j=0}^{7}
P^{(k-1)}_{8n+j},
\qquad
n\in\mathbb{N}.
\label{pattern:eq1.1}
\end{equation}

We define the indicator function

\begin{equation}
\phi(S)=
\begin{cases}
1 & \text{if } S \text{ contains at least one prime},\\
0 & \text{otherwise}.
\end{cases}
\label{eq:phiMoins2}
\end{equation}

We extend this definition to encoded integers by setting

\begin{equation}
\phi(n)=
\begin{cases}
1 & n \neq 0,\\
0 & \text{otherwise}.
\end{cases}
\label{eq:phiMoins3}
\end{equation}

The set $P^{(1)}_n$ is encoded by the integer in $[0, 255]$ whose binary representation is
\[
f^{(1)}(P^{(1)}_n)
=
(\phi(P^{(0)}_{8n+1})\cdots
\phi(P^{(0)}_{8n+8}))_2.
\]

More generally,

\[
f^{(k)}(P^{(k)}_n)
=
(\phi(P^{(k-1)}_{8n})\cdots
\phi(P^{(k-1)}_{8n+7}))_2.
\]

If
$f^{(k)}(P^{(k)}_n)\neq0$,
then at least one of the eight sub-blocks
$P^{(k-1)}_{8n+j}$,
$j=0,\ldots,7$,
contains at least one prime number.
If
$f^{(k)}(P^{(k)}_n)=255$,
then every sub-block contains at least one prime number.

Throughout the paper, the binary encoding $f^{(k)}(P^{(k)}_n)$ is called the \emph{$k$-pattern}.

For notational convenience, we write $f^{(k)}(n)$ instead of
$f^{(k)}(P^{(k)}_n)$, and $\phi(f^{(k)}(n))$ instead of $\phi(P^{(k)}_n)$.

Let us detail the first three levels of the multi-scale representation :
	
\paragraph{Level 1 (fine scale)}The set $\mathbb{N}^*$ is partitioned into blocks of 8 consecutive integers $P^{(1)}(n) = \{8n+1,\ldots,8n+8\}$. Each block $P^{(1)}(n)$ is encoded in binary by $\phi(\{8n+1\})\cdots \phi(\{8n+8\})$ (1 if $8n+j$ is prime and 0 otherwise), then converted to decimal. For example: $(1,2,\ldots,8) \rightarrow (01101010)_2 = 106$ and $(9,\ldots,16) \rightarrow (00101000)_2 = 40$.\\	
The first terms of this sequence are: $106, 40, 162, 10, 8, 162, 8, 40, 34, 130, \ldots$

\paragraph{Level 2 (intermediate scale)}
At this level $\mathbb{N}^*$ is partitioned with the sets  $P^{(2)}(n)$. We aggregate the terms $P^{(1)}(8n+j), j=0,\ldots,7$,  from level 1 into blocks of 8.  We encode $P^{(2)}(n)$ by
	$f^{(2)}(n)=f^{(2)}(P^{(2)}(n))= (\phi(P^{(1)}_{8n})\ldots\phi(P^{(1)}_{8n+7}))_2 = 
	(\phi(f^{(1)}(8n))\cdots \phi(f^{(1)}(8n+7)))_2$
	
	 For example,
	$(f^{(1)}(0),f^{(1)}(1),\ldots,f^{(1)}(7)) =(106,40, 162,\ldots,8,40) \rightarrow 
	f^{(2)}(0)=(\phi(f^{(1)}(0))\cdots \phi(f^{(1)}(7)))_2 = (11111111)_2 = 255$. $f^{(2)}(0)$ encodes 64 integers, it shows that each sub-block of 8 integers contains at least one prime.
A level-2 pattern such as $95 = (01011111)_2$ indicates that 6 of the 8 sub-blocks contain at least one prime.
	
\begin{figure}[h]
		\centering
		\resizebox{1.0\linewidth}{!}{\begin{tikzpicture}[x=1.5cm, y=1cm]
  \node[] (racine) at (0,0) {255};
  
  \foreach \i [count=\x from 0] in {255,255,255,191,251,127,255,255} {
    \node[] (n1-\x) at (-7+\x*2, -2) {\i};
    \draw (racine) -- (n1-\x);
  }
  
  \foreach \i [count=\x from 0] in {138,0,32,2,40,130,128,32} {
    \node[] (n2-\x) at (-6.5+\x*1.5, -4) {\i};
    \draw (n1-3) -- (n2-\x);  
  }
  
  \foreach \i [count=\x from 0] in {225,226,227,228,229,230,231,232} {
    \node[] (n3-\x) at (-5.5+\x*1, -6) {\i};
    \draw (n2-4) -- (n3-\x);  
  }
\end{tikzpicture}}
		\caption{Multi-scale tree illustrating three resolution levels. Each node encodes the presence of prime numbers in a specific interval: 512 integers (level 3), 64 integers (level 2) and 8 integers (level 1). The value $40$ at level 1 corresponds to integers $225-232$, of which only $227$ and $229$ are prime.}
		\label{fig:arbre}
	\end{figure}
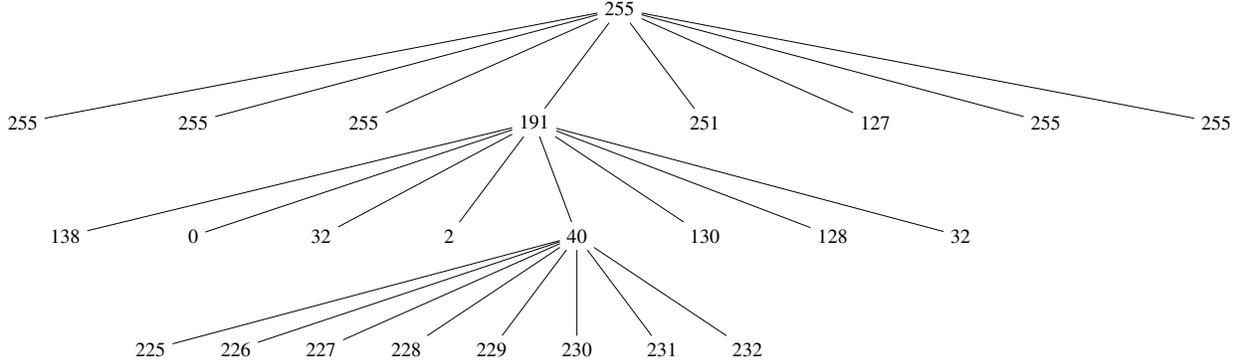

\paragraph{Level 3 (coarse scale)}
	By iterating the process at level 2, we obtain a coarser-scale view. Each term of the sequence encodes primality in a block of $8^3=512$ integers. The first 366 terms are equal to $255$, indicating that each of the corresponding sub-blocks contains at least one prime. As shown in Figure~\ref{fig:arbre}, the hierarchy takes the form of a tree, with information growing increasingly refined at lower levels:
	
	\begin{itemize}
		\item The root node $255$ indicates that the first 512 integers contain at least one prime in each 64-integer sub-block.
		 
		\item The node $191$ at level 2 already reveals a slight irregularity
		\item The node $40$ at level 1 allows the primes $227$ and $229$ to be localized within the corresponding 8-integer block.
	\end{itemize}

	Let us now formally describe this multi-scale representation. The hierarchical construction is performed via the sequences $(f^{(k)})_{k\geq 1}$ defined recursivelyt by:
	
	\begin{equation}
		\begin{aligned}
			f^{(1)}(n) &= \sum_{j=0}^{7} \phi(\{8n+1+j\}) \cdot 2^{7-j}, \quad n \in \mathbb{N}\\
			f^{(k+1)}(n) &= \sum_{j=0}^{7} \phi \big( f^{(k)} (8n+j) \big) \cdot 2^{7-j}, \quad n \in \mathbb{N}
		\end{aligned}
		\label{Scale:eq1}
	\end{equation}

Based on this multi-scale representation, we investigate the following question. 
Given an integer $m$, we study the quantity $\pi_m^{(k)}(n)$, for $n=0,\ldots,255$, defined by

\[
\pi_m^{(k)}(n)
=
\operatorname{card}
\left\{
P^{(k)}_j \mid
(\phi(P^{(k-1)}_{8j})\phi(P^{(k-1)}_{8j+1})\cdots
\phi(P^{(k-1)}_{8j+7}))_2=n,\;
j=0,\ldots,\left\lfloor\frac{m}{8^k}\right\rfloor
\right\}.
\]

Equivalently,
\[
\pi_m^{(k)}(n)
=
\operatorname{card}
\left\{
P^{(k)}_j \mid
f^{(k)}(j)=n,\;
j=0,\ldots,\left\lfloor\frac{m}{8^k}\right\rfloor
\right\}.
\]

The proposed multi-scale representation provides a natural way to define the quantities $\pi_m^{(k)}(n)$. The main objective of this article is to derive estimates for these quantities (see Observations~2).

Furthermore, for sufficiently large values of $m$, we observe that the histograms of $\{\pi_m^{(k)}(n)\}_{n=0}^{255}$ consistently exhibit the same invariant shape. This persistent profile suggests a structural feature of the distribution of prime numbers.
	
\section{The sequences generated by the multi-scale representation}
\subsection{The sequence at the first level}
The sequence $(f^{(1)}(n))$ generated at the first level of the multi-scale representation, see (\ref{Scale:eq1}), encodes primality in blocks of 8 consecutive integers:
$f^{(1)}(n)$ is an integer between 0 and 255. Its binary representation defines a pattern (1-pattern) that identifies which of the 8 consecutive integers are prime and which are composite.

 The first 49 values of $f^{(1)}(n)$ are organized in a spiral, see Figures \ref{spiral} and \ref{F1Premier}. A quick look at the sequence makes it easy to visualize the distribution of prime numbers. For example, 40 indicates the presence of a twin prime, while the other numbers in the corresponding block are composite. We first make the following observation about the sequence.
	
\paragraph{Pattern Restrictions.}
\paragraph{Observation 1.} \textit { The integer $f^{(1)}(n)$ takes only 14 distinct values:
\begin{equation}
f^{(1)}(n) \in C
\end{equation}
\centerline{where $C = \{0,2,8,10,32,34,40,106,128,130,136,138,160,162\}$ }}
	
\noindent \textit{Proof.} The absence of certain patterns is explained by arithmetic constraints. For example, a pattern like $(00101010)_2 = 42$ is impossible because among three consecutive odd numbers, at least one is divisible by 3.\\

In the remainder of the article, we denote by $np(j)$ the function that counts the number of ones in the binary representation of $j$.

\paragraph{Distribution of primes per block.} The number of primes within a 1-pattern is given by:
	
	\begin{equation}
		np(f^{(1)}(n)) = \begin{cases}
			0 & \text{if } f^{(1)}(n) = 0 \\
			1 & \text{if } f^{(1)}(n) \in \{2,8,32,128\} \\
			2 & \text{if } f^{(1)}(n) \in \{10,34,40,130,136,160\} \\
			3 & \text{if } f^{(1)}(n) \in \{138,162\} \\
			4 & \text{if } f^{(1)}(n) = 106
		\end{cases}
	\end{equation}

\paragraph{Interpretation of the elements of the sequence $f^{(1)}(n)$.}
The sequence $f^{(1)}(n)$ provides an encoding for a classification of blocks according to their prime-number patterns into subsets $E_t$, where $t \in C$.  
Each subset is defined by $\{n\in\mathbb N \mid f^{(1)}(n)=t\}$.
Here are some examples of these subsets $E_t$:

\begin{itemize}
\item $E_2=\{n\in\mathbb N \mid 8n+7 \text{ is prime, while } 8n+1,8n+3,\text{ and }8n+5\text{ are not prime}\}$
\item $E_{128} = \{n \in \mathbb{N} \mid 8n+1 \text{ is prime, while } 8n+3, 8n+5, \text{ and } 8n+7 \text{ are not prime}\}$
\item $E_{10} = \{n \in \mathbb{N} \mid 8n+1 \text{ and } 8n+7 \text{ are prime, while } 8n+3 \text{ and } 8n+5 \text{ are not prime}\}$
\item $E_{138} = \{n \in \mathbb{N} \mid 8n+1, 8n+5 \text{ and } 8n+7 \text{ are prime, but } 8n+3 \text{ is not prime}\}$
\item Similarly, the subsets $E_t$ can be defined for $t \in \{0,8,32,34,40,130,136,160,162\}$
\end{itemize}
	
	\begin{figure}[h]
		\centering
		\begin{minipage}{0.4\textwidth}
			\centering
			\resizebox{0.7\linewidth}{!}{\begin{tikzpicture}[scale=1.5]
  \draw[gray!30] (0,0) grid (7,7);
  
  
  \draw[red, thick, ->] (3.5,3.5) -- (4.5,3.5) -- (4.5,4.5) -- (3.5,4.5) -- 
                         (2.5,4.5) -- (2.5,3.5) -- (2.5,2.5) -- (3.5,2.5) -- 
                         (4.5,2.5) -- (5.5,2.5) -- (5.5,3.5) -- (5.5,4.5) -- 
                         (5.5,5.5) -- (4.5,5.5) -- (3.5,5.5) -- (2.5,5.5) -- 
                         (1.5,5.5) -- (1.5,4.5) -- (1.5,3.5) -- (1.5,2.5) -- 
                         (1.5,1.5) -- (2.5,1.5) -- (3.5,1.5) -- (4.5,1.5) -- 
                         (5.5,1.5) -- (6.5,1.5);
  
  \fill[blue] (3.5,3.5) circle (3pt);
  \node[above right] at (3.5,3.5) {};
\end{tikzpicture}}
			\caption{Spiral arrangement of sequence elements}
			\label{spiral}
		\end{minipage}
		\hspace {0.09 cm}
		\begin{minipage}{0.45\textwidth}
			\centering
			\begin{tabular}{*{7}{>{\centering\arraybackslash}m{0.5cm}}}
				34 &   0 &   8 &  10 & 130 &  32 & 128 \\ 
				136 &  32 &   2 & 128 &  40 & 138 & 130 \\ 
				0 & 160 &   8 &  10 & 162 & 128 &  40 \\ 
				32 &  10 & 162 & \textbf{106} &  40 &  32 &   2 \\ 
				0 &   8 &   8 &  40 &  34 & 130 &  32 \\ 
				32 &  34 &   8 &  40 &   2 & 138 &   0  \\ 
				128 &  40 & 130 &   2 &   8 &  34 &   8 \\ 
			\end{tabular}
			\caption{First 49 elements of the sequence $f^{(1)}(n)$ arranged in a spiral.}
			\label{F1Premier}
		\end{minipage}
		
		\vspace{0.5cm} 
		\begin{minipage}{0.45\textwidth}
			\centering
			\begin{tabular}{*{7}{>{\centering\arraybackslash}m{0.3cm}}}
				158& 223& 157& 182& 239& 254& 83 \\ 
				246&  221& 251& 111& 182& 183& 239 \\ 
				220& 253& 251& 191& 255& 255& 122 \\ 
				159& 255& 127& \textbf{255}& 255& 239& 243 \\ 
				59& 237& 255& 255&  95& 254& 157 \\ 
				111& 252&  58& 255& 123& 253& 249  \\ 
				247& 222& 189& 189& 115& 223& 238 \\ 
			\end{tabular}
			\caption{First 49 elements of the sequence $f^{(2)}(n)$ arranged in a spiral. }
			\label{F2Premier}
		\end{minipage}
		\hspace{0.10 cm} 
		\begin{minipage}{0.45\textwidth}
			\centering
			\begin{tabular}{*{7}{>{\centering\arraybackslash}m{0.3cm}}}
				255 &  255 &  255 & 255 & 255 &  255 & 255\\ 
				255 & 255 & 255 &  255&   255 &  255 &  255 \\ 
				255 & 255 & 255 &  255&   255 &  255 &  255 \\ 
				255& 255 &  255 & 255 &  255 & 255 & 255 \\ 
				255 & 255 & 255 &  255&   255 &  255 &  255 \\ 
				255 & 255 & 255 &  255&   255 &  255 &  255  \\ 
				255 & 255 & 255 &  255&   255 &  255 &  255  \\ 
			\end{tabular}
			\caption{$f^{(3)}(n)$ defined in Example 1 is constant up to $n = 366$}
		\end{minipage}
	\end{figure}
In the remainder of the article, the study is conducted using the files containing the prime numbers less than $982\;451\;653$ available on the website
\href{https://t5k.org/lists/small/millions/}{{\em The first fifty million primes}} originally created by Chris Caldwell \cite{Caldwell, Caldwell2009}.
Our study is done over the interval $[1, N_{\max}]$. We choose $N_{\max}=982\,451\,200$, the largest multiple of $8^3$ less than $982\,451\,653$, 
so that the interval $[1,N_{\max}]$ can be exactly partitioned into blocks of sizes $8$, $8^2$, and $8^3$, corresponding to the first three levels of the multi-scale representation.

\begin{figure}[h]
\centering
\begin{minipage}[t]{0.45\textwidth}
\centering
    \begin{tabular}{lc@{\hspace{2em}}!{\vrule width 0.4pt}@{\hspace{0.5em}}!{\vrule width 0.4pt}@{\hspace{1em}}lc}
        \toprule
        $n$ & $\pi^{(1)}_{N_{\max}}(n)$ & $n$ & $\pi^{(1)}_{N_{\max}}(n)$ \\
        \midrule
        0   & 78\,517\,574    & 106 & 1            \\
        2   & 9\,316\,943     & 128 & 9\,314\,958     \\
        8   & 10\,066\,986    & 130 & 1\,497\,811     \\
        10  & 748\,618       & 136 & 749\,060      \\
        32  & 10\,065\,624    & 138 & 93\,887      \\
        34  & 749\,468       & 160 & 750\,314        \\
        40  & 842\,208       & 162 & 92\,948       \\
        \bottomrule
    \end{tabular}
    \captionof{table}{Prime numbers: Histogram $\pi^{(1)}_{N_{\max}}(n)$ of 1-patterns}
    \label{tab1:distribution}
\end{minipage}
\vspace{1 cm}
\begin{minipage}[t]{0.45\textwidth}
			\begin{tabular}{lc}
				\toprule
				Configuration Type & Proportion \\
				\midrule
				Blocks without primes & 63.94\% \\
				Blocks with 1 prime & 31.56\% \\
				Blocks with 2 primes & 4.35\% \\
				Blocks with 3-4 primes & 0.15\% \\
				\bottomrule
			\end{tabular}
	\captionof{table}{Prime numbers: Distribution of 1-patterns}
	\label{tab2:distribution}
		\end{minipage}
\end{figure}
\paragraph{Counting 1-patterns.} The statistical distribution of the sequence $f^{(1)}(n)$ is presented in Tables \ref{tab1:distribution} and \ref{tab2:distribution}. Over the interval $[1, m]$, we observe that 

the number of empty 1-patterns increases more rapidly than the number of 1-patterns containing a single prime, while the latter increases more slowly. 
Moreover, the number of 1-patterns containing two or three prime numbers decreases as $m$ increases. Table \ref{tab2:distribution} shows this distribution for $m = N_{\max}$. 

\subsection{The sequence at the second level}
The sequence $f^{(2)}(n)$ takes all values between 0 and 255. The first observed values are $255, 255, 255, 191, 251, 127,$ and $255$. The first 49 values of this sequence are organized in a spiral and represented in Figure~\ref{F2Premier}. 

\begin{figure}[h]
    \centering
    \begin{subfigure}{0.48\textwidth}
        \centering
        \includegraphics[width=\linewidth]{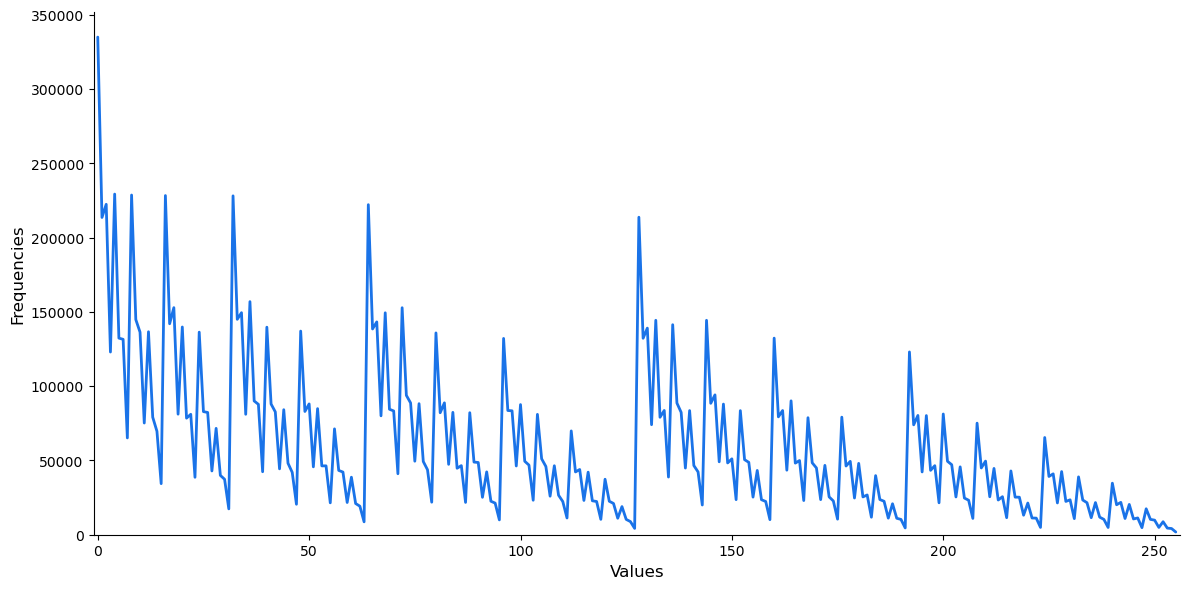}
        \caption{Histogram of $f^{(2)}$ values}
        \label{fig:f2_histogram}
    \end{subfigure}
    \hfill
    \begin{subfigure}{0.48\textwidth}
        \centering
        \includegraphics[width=\linewidth]{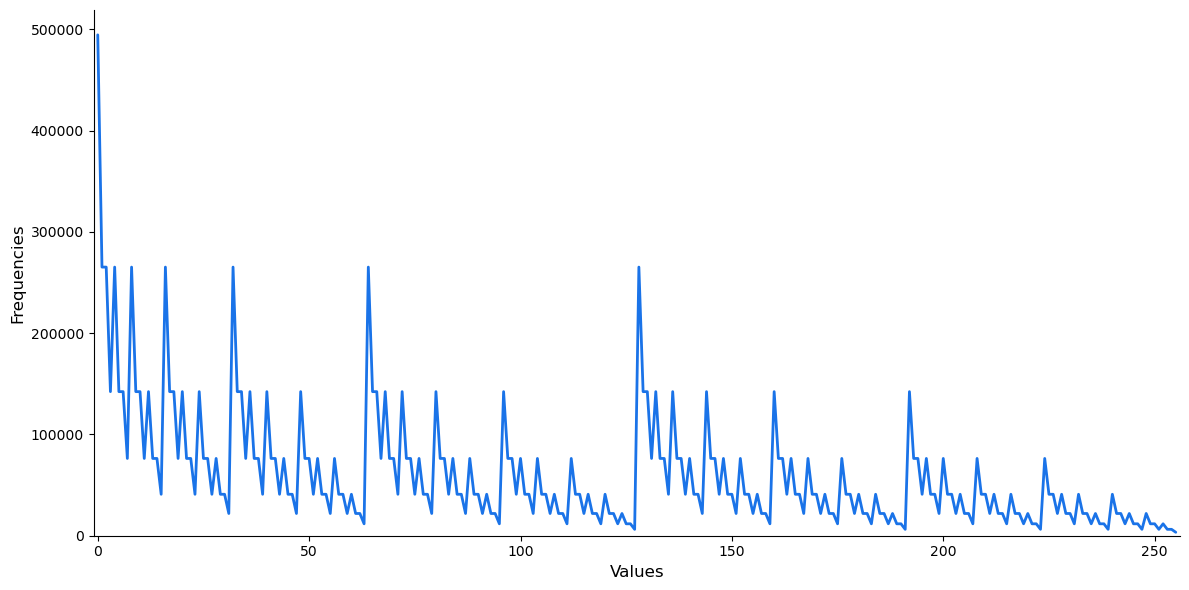}
        \caption{Theoretical histogram with $C^{(2)}_m = 1$ and $\epsilon^{(2)}_j(m) = 0$ (Cf. observation 2)}
        \label{fig:theoretical_histogram}
    \end{subfigure}
    \caption{Shapes of empirical and theoretical histograms for $f^{(2)}$, computed over $[1,N_{\max}]$.}
    \label{fig:histograms}
\end{figure}

\paragraph{The Invariant Shape of $k$-Pattern Histograms} 
The distribution of the sequence $f^{(2)}(n)$ for the integers $n$ in the interval $[1,\frac{N_{\max}}{8^2}]$ is shown in Figure~\ref{fig:f2_histogram}, 
whereas Figure~\ref{fig:theoretical_histogram} illustrates the histogram obtained using the model given by Observation~2. 
The histogram reveals an irregular but repeating shape. 

The observed periodicity stems from the statistical predominance of 2-patterns with $np(f^{(2)}(n))=1$ over those with $np(f^{(2)}(n))=2$, and similarly for higher values.
The characteristic shape of the histogram becomes apparent for $m>16\,000\,000$. We conjecture that, for a sufficiently large interval, the shape illustrated in Figure~\ref{fig:histograms} is observed for the histograms of the sequences $f^{(k)}(n)$ with $k\geq2$. Moreover, this shape remains invariant as the interval size increases. This conjecture is stated as Observation~2.

\subsection{The sequence at the third level}
	\label{subsec:niveau3}
	
	\paragraph{Pattern restrictions.} 
	Over the interval $[1, N_{\max}]$, only 114 distinct values are observed among the 256 possible values:

	\[ 
	\begin{aligned}
		f^{(3)}(n) \in
 \{ &31, 45, 47, 51, 54, 55, 58, 59, 61, 62, 63, 70, 71, 77, 78, 79, 83, 85, 86, 87, \\
 &89, 90, 91, 93, 94, 95, 103, 107, 109, 110, 111, 115, 117, 118, 119, 121,\ldots, 127, \\
 &142, 143, 151, 155, \ldots, 159, 165, 167, 169, \ldots, 175, 179, \ldots, 183, 185, \ldots, 191, 197,199, \\
 &203, 205, 206, 207, 211, 213, 214, 215, 217, \ldots, 223, 226, 227, 229,230, 231, 233, \ldots, 239, 241, \ldots, 255\}
	\end{aligned}
	\]
	
The histogram of the 3-patterns $f^{(3)}(n)$ for $n\in[0,N_{\max}/512]$ is given in Appendix Table~\ref{tab:hp3}. We observe that the value 255 represents 83.78\% of the 3-patterns. All these 3-patterns contain at least eight prime numbers. 

The sequence $f^{(3)}(n)$ remains equal to $255$ for all $n\leq366$, with the first deviation occurring at $f^{(3)}(367)=223$.

\paragraph{Estimates for the number of $k$-patterns}
Below, we propose an expressions for estimating the quantities $\pi^{(k)}_{m}(n)$:

\paragraph{Observation 2.} {\em

\begin{itemize}

\item [(1)] The probability that an odd integer less than or equal to $m$ is prime is approximately $q_0(m) = \frac{2 \text{li}(m)}{m}$, where
\begin{equation*}
    \operatorname{li}(x) = \int_{2}^{x} \frac{1}{\ln(t)} dt
    \label{eq:logarithmic_integral}
\end{equation*} 

\item [(2)] The probability that a $1$-pattern contains at least one prime is approximately equal to

\[ 
    q_1(m) = 1 - (1 - q_0(m))^4  
\]

\item [(3)] For $k>1$, the probability that a $k$-pattern contains at least one prime number is approximately equal to

    \[ 
    q_k(m) = 1 - (1 - q_{k-1}(m))^8  =  1 - \big(1-q_0(m)\big)^{4 \times 8^{k-1}}
    \]

\item [(4)] Under the assumption of independence, the number $\pi^{(k)}_{m}(n)$ of $k$-patterns with a given configuration $n$ is approximately

    \[
    \pi^{(k)}_{m}(n) = \frac{m}{8^k} q_{k-1}(m)^j \big(1-q_{k-1}(m)\big)^{8-j}
    \]
    \centerline{where $np(n) = j$, }

The numerical results show systematic deviations from this approximation. To account for these deviations, we conjecture that

    \[
    \pi^{(k)}_{m}(n) \approx ( C^{(k)}_m + \epsilon^{(k)}_j(m) )\frac{m}{8^k} q_{k-1}(m)^j \big(1-q_{k-1}(m)\big)^{8-j}
    \]
    \centerline{where $np(n) = j$, }
$C^{(k)}_m$ is a constant, and $\epsilon^{(k)}_j(m)$ represents a fluctuation around $C^{(k)}_m$, satisfying $| \epsilon^{(k)}_j(m) | \ll C^{(k)}_m $.

\end{itemize}
}

Before providing the justification for these observations, we first report some numerical results for $m=N_{\max}$: 
We obtain $q_0(N_{\max})=0.101790$, while the observed value is $q_0^{\mathrm{obs}}(N_{\max})=0.101783$, corresponding to an absolute error of approximately $10^{-6}$.

For $q_1$, we obtain $q_1(N_{\max})=0.349105$, while the observed value is $q_1^{\mathrm{obs}}(N_{\max})=0.360639$, corresponding to an absolute error of approximately $10^{-2}$.

For $q_2$, we obtain $q_2((N_{max}) =0.967783$, while the observed value is $q_2^{\mathrm{obs}}(N_{\max})=0.978178$, corresponding to an absolute error of approximately $10^{-2}$.

For $q_3$, we obtain $q_3((N_{max}) =0.999999$, while the observed value is $q_3^{\mathrm{obs}}(N_{\max})=1$. Although $q_3(N_{\max})$ is displayed as $0.999999$, its difference from the observed value is approximately $10^{-12}$ before rounding.

For $m=N_{\max}$, levels beyond $k=3$ provide no additional information within the accuracy considered.

\noindent {\em Justification.} 
\begin{itemize}
\item [(1)]
For the odd integers not exceeding $m$, the proportion that are prime is given by
$\frac{\pi(m)}{m/2} = \frac{2\pi(m)}{m}$,
we approximate this quantity by $\frac{2 \operatorname{li}(m)}{m}$. Where $\pi(m)$ denotes the number of primes less than or equal to $m$.

\item [(2)] The probability that a 1-pattern contains no prime is 
$(1 -q_0(m))^4$, because a 1-pattern contains four odd positions that can contain prime numbers. Therefore, the probability that a 1-pattern contains at least one prime is 
$q_1(m) = 1-(1 - q_0(m))^4$.

\item [(3)] Consider the block $P^{(k)}_n=\bigcup_{j=0}^7P^{(k-1)}_{8n+j}$ associated with a $k$-pattern.
The probability that none of the blocks $P^{(k-1)}_{8n+j}$, $j=0,\ldots,7$, contains a prime is $(1-q_{k-1}(m))^8$. By recurrence, we can write $(1-q_{k-1}(m))^8 = (1-q_0(m))^{4\times 8^{k-1}}$.
Therefore, the probability that at least one of the blocks $P^{(k-1)}_{8n+j}$ contains at least one prime number is \[1-(1-q_{k-1}(m))^8 \]

\item [(4)] To establish the formulas in item (4), 
we illustrate the argument using $\pi_m^{(1)}(j)$ for $j\in\{2,8,32,128\}$.
The same reasoning applies to $\pi^{(k)}_{m}(j)$, with $k \geq 2$. Consider the case $j = 2$ where the 1-pattern $f^{(1)}(n) = j$ corresponds to $8n + 7$ being prime while $8n+1$, $8n+3$ and $8n+5$ are not, for $n \leq \frac{m}{8}$. Assuming independence, the probability that $8n+7$ is the only prime number in the sequence is:
\[ q_0(m)(1-q_0(m))^3 \]

This formula remains valid for $j \in \{2,8,32,128\}$.

In the interval $[1, m]$, the number of 1-patterns is approximately equal to $\frac{m}{8}$, so the count $\pi^{(1)}_{m}(j)$ can be written as:

\[ \pi^{(1)}_{m}(j) \approx \frac{m}{8} q_0(m)(1-q_0(m))^3 \]

To correct the approximations and the independence assumption, we propose a correction as follows:
\[ \pi^{(1)}_{m}(j) = (C^{(1)}_m + \epsilon^{(1)}_j(m)) \frac{m}{8} q_0(m)(1-q_0(m))^3 \]
where $\epsilon^{(1)}_j(m)$ represents a fluctuation around the constant $C^{(1)}_m$.

\end{itemize}

	\section{Conclusion}
	\label{sec:conclusion}
	This work explores a multi-scale approach to studying arithmetic properties of natural numbers. Through the example of prime numbers, we illustrate how this hierarchical representation can reveal structural features in the distribution of prime numbers. The main features of this method include encoding arithmetic properties into nested sequences, enabling the simultaneous analysis of local and global behaviors.

In this article, the generated sequences take integer values between 0 and 255. This bounded representation provides a compact encoding of prime number distributions across intervals of large integers. We further analyze these distributions through scale-dependent histograms, each comprising at most 256 bins. Each bin counts the occurrences of a specific pattern encoding the presence or absence of primes within the corresponding blocks. Based on a probabilistic model, we propose an expression for estimating the bin counts in these histograms. For sufficiently large intervals, empirical results suggest that the shape of these histograms stabilizes and approaches an apparently invariant form. This persistent structure suggests the existence of a characteristic pattern in the multi-scale representation of prime number distributions.

Future research directions may include refining the proposed probabilistic estimates for histogram bins, exploring multi-scale representations with overlapping pattern schemes, and adapting the method to variable pattern sizes.


\section{Appendix}

\begin{table}[h]
\begin{tabular}{|c|c||c|c||c|c||c|c||c|c||c|c|}
\hline
$n$ & $\pi^{(3)}_{N_{\max}}(n)$ & $n$ & $\pi^{(3)}_{N_{\max}}(n)$ & $n$ & $\pi^{(3)}_{N_{\max}}(n)$ & $n$ & $\pi^{(3)}_{N_{\max}}(n)$ & $n$ & $\pi^{(3)}_{N_{\max}}(n)$ & $n$ & $\pi^{(3)}_{N_{\max}}(n)$ \\
\hline
31 & 15 & 45 & 1 & 47 & 17 & 51 & 1 & 54 & 1 & 55 & 20 \\
58 & 2 & 59 & 17 & 61 & 15 & 62 & 10 & 63 & 649 & 70 & 1 \\
71 & 1 & 77 & 4 & 78 & 2 & 79 & 8 & 83 & 3 & 85 & 1 \\
86 & 1 & 87 & 17 & 89 & 1 & 90 & 1 & 91 & 13 & 93 & 18 \\
94 & 17 & 95 & 830 & 103 & 18 & 107 & 23 & 109 & 13 & 110 & 20 \\
111 & 761 & 115 & 20 & 117 & 26 & 118 & 15 & 119 & 798 & 121 & 15 \\
122 & 14 & 123 & 820 & 124 & 14 & 125 & 837 & 126 & 868 & 127 & 35878 \\
142 & 1 & 143 & 13 & 151 & 9 & 155 & 12 & 156 & 1 & 157 & 16 \\
158 & 20 & 159 & 709 & 165 & 1 & 167 & 17 & 169 & 1 & 170 & 1 \\
171 & 15 & 172 & 1 & 173 & 18 & 174 & 22 & 175 & 788 & 179 & 15 \\
181 & 24 & 182 & 19 & 183 & 804 & 185 & 11 & 186 & 22 & 187 & 824 \\
188 & 19 & 189 & 824 & 190 & 798 & 191 & 35995 & 197 & 1 & 199 & 10 \\
203 & 14 & 205 & 20 & 206 & 21 & 207 & 704 & 211 & 21 & 213 & 12 \\
214 & 21 & 215 & 763 & 217 & 14 & 218 & 16 & 219 & 832 & 220 & 14 \\
221 & 792 & 222 & 835 & 223 & 36046 & 226 & 1 & 227 & 13 & 229 & 16 \\
230 & 12 & 231 & 715 & 233 & 20 & 234 & 22 & 235 & 832 & 236 & 15 \\
237 & 845 & 238 & 835 & 239 & 35953 & 241 & 14 & 242 & 9 & 243 & 699 \\
244 & 23 & 245 & 771 & 246 & 775 & 247 & 35991 & 248 & 8 & 249 & 708 \\
250 & 784 & 251 & 36169 & 252 & 671 & 253 & 36182 & 254 & 36171 & 255 & 1607654 \\
\hline
\end{tabular}
\caption{Histogram $\pi^{(3)}_{N_{\max}}(n)$ of 3-patterns}
\label{tab:hp3}
\end{table}


\begin{thebibliography}{1}\footnotesize
	 	
	 	\bibitem{Iovane2009} Gerardo Iovane, The set of prime numbers: Multiscale analysis and numeric accelerators, {\it Journal Chaos, Solitons \& Fractals, Elsvier}, {\bf 41}(4) (2009), 1953-1965, 
	 	
	 	\bibitem{dusart2018} P. Dusart, Explicit estimates of some functions over primes, {\it The Ramanujan Journal},{\bf 45}(1) (2018), 227-251. 
	 	
	 	\bibitem{goldston2009}	D. A. Goldston, J. Pintz,  and C. Y. Y{\i}ld{\i}r{\i}m, Primes in Tuples I, {\it Annals of Mathematics}, {\bf 170}(2) (2009), 819-862. 
	 	
	 	\bibitem{maynard2015}	J. Maynard, Small gaps between primes, {\it Annals of Mathematics}, {\bf 181}(1) (2015), 383-413.
	 	
	 	\bibitem{Caldwell}\url{https://en.wikipedia.org/wiki/PrimePages}
 	
	 	\bibitem{Caldwell2009} Chris K. Caldwell and Jr., G. L. Honaker, Prime Curios!: The Dictionary of Prime Number Trivia, {\it CreateSpace Independent Publishing Platform}, 2009
 	
	 	\bibitem{Tenenbaum}	G. Tenenbaum, Introduction to Analytic and Probabilistic Number Theory, {\it Cambridge Studies in Advanced Mathematics}, 1995. 
		
		
	\end{thebibliography}
\end{document}